\documentclass[12pt]{article}
\usepackage[utf8]{inputenc}
\usepackage{amssymb,amsfonts}
\usepackage[english,russian]{babel}
\usepackage[usenames]{color}
\usepackage{wrapfig}
\usepackage{url}
\usepackage{amsmath}

\newcommand{\3}{''}
\newcommand{\q}{\mathbb{Q}}

\newcommand{\n}{\mathbb{N}}
\newcommand{\R}{\mathbb{R}}
\newcommand{\fo}[1]{${\displaystyle #1}$}
\newcommand{\hor}[1]{\q\left[\sqrt{#1}\right]}
\newcommand{\horp}[1]{\q^+\left[\sqrt{#1}\right]}
\newcommand{\x}{\overline{x}}
\newcounter{teorema}
\newcounter{glav}
\newcounter{risunok}
\newcommand{\ris}{\refstepcounter{risunok}{\bf Figure \arabic{risunok}}}
\newcommand{\glava}{\refstepcounter{glav}\begin{center}{\bf\arabic{glav}.}}
\newcommand{\teo}{\refstepcounter{teorema}{\bf Theorem \arabic{teorema}}}
\newcounter{le}
\newcounter{op}
\newcommand{\lemma}{\refstepcounter{le}\textbf{Lemma \arabic{le}. }}
\newcommand{\opred}{\addtocounter{op}{1}\textbf{Definition \arabic{op}. }}
\begin{document}
\begin{center}
{\Large \textbf{Dissection of a rectangle into rectangles with given side ratios}}\\
$\phantom{kl}$\\
\begin{large}
\textbf{Fyodor Sharov\footnote{The author was partially supported by the President of the Russian Federation grant MK-6137.2016.1.}}\\
$\phantom{kl}$\\
\end{large}
{\small\textit{Faculty of Mathematics, National Research University Higher School of Economics}}
\end{center}
{\small 
\begin{center}
\textbf{Abstract}
\end{center}
\renewcommand{\refname}{References}

Given a collection of $n$ rectangles such that the side ratio of each one is a quadratic irrationality, we find all rectangles which can be tiled by rectangles similar to one of the given ones. It means that each possible shape can be used several times or not used at all, so that the number of rectangles in the tiling is not necessarily equal to $n$.


\glava\textbf{ Main result.}\end{center}

The aim of this paper is to prove the following result.

\teo\textbf{.} \textit{Consider numbers ${x_1=a_1+b_1\sqrt{p}}$,~\dots,~${x_n=a_n+b_n\sqrt{p}}$ such that for any \fo{i\in\lbrace 1,2, \dots, n\rbrace} we have ${x_i>0}$, $a_i,b_i,p\in\q$ and ${\sqrt{p}\notin\q}$.\\
1) If there exist \fo{i,j\in\lbrace 1,2, \dots, n\rbrace} such that ${(a_i-b_i\sqrt{p})(a_j-b_j\sqrt{p})<0}$, then a rectangle with side ratio $z$ can be tiled by rectangles with side ratios $x_1$,~\dots,~$x_n$ if and only if $${\displaystyle z\in\left\lbrace  e+f\sqrt{p}>0\,\vert\, e,f\in\q \right\rbrace.}$$}
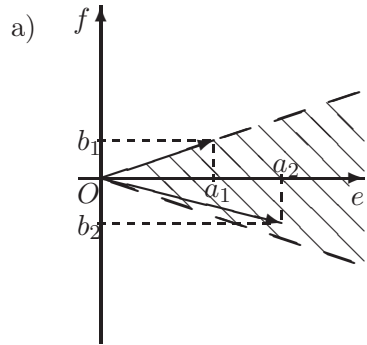
\label{mainT}
\begin{wrapfigure}[20]{l}{0.27\linewidth} 
\vspace{-10ex}
\hspace*{2.75ex}
\begin{picture}(45,145) 
\thicklines
\put(-5,141){a)}
\put(-5,86){b)}
\put(0,55){\vector(1,0){42}}
\put(22,52){\vector(0,1){38}}
\put(22,55){\vector(1,3){5}}
\put(22,55){\vector(-1,4){6}}
\multiput(27.3,70.9)(2.4,7.2){3}{\line(1,3){1.2}}
\multiput(21.7,55.9)(-2.4,7.2){5}{\line(-1,3){1.2}}
\multiput(27,69.5)(0,-2){8}{\line(0,-1){1}}
\multiput(26.5,70)(-2,0){3}{\line(-1,0){1}}
\multiput(16,77.5)(0,-2){12}{\line(0,-1){1}}
\multiput(16.5,79)(2.4,0){3}{\line(1,0){1.2}}
\thinlines
\put(20,61){\line(1,-1){3}}
\put(18,67){\line(1,-1){6}}
\put(16,73){\line(1,-1){9}}
\put(14,79){\line(1,-1){12}}
\put(12,85){\line(1,-1){15}}
\put(11,90){\line(1,-1){17}}
\put(15,90){\line(1,-1){14}}
\put(21,88){\line(1,-1){9}}
\put(23,90){\line(1,-1){8}}
\put(27,90){\line(1,-1){5}}
\put(31,90){\line(1,-1){2}}
\put(25.8,52){$a_1$}
\put(40.4,51.4){$e$}
\put(14.8,52){$a_2$}
\put(18.8,51.7){$O$}
\put(18.9,68.3){$b_1$}
\put(22.5,77.9){$b_2$}
\put(18.4,87.6){$f$}

\put(-48,100){\begin{picture}(100,45)

\thicklines
\put(55,0){\vector(0,1){45}}
\put(52,22){\vector(1,0){38}}
\put(55,22){\vector(3,1){15}}
\put(55,22){\vector(4,-1){24}}
\multiput(70.9,27.3)(7.2,2.4){3}{\line(3,1){4}}
\multiput(55.9,21.7)(7.2,-2.4){5}{\line(3,-1){4}}
\multiput(69.5,27)(-2,0){8}{\line(-1,0){1}}
\multiput(70,26.5)(0,-2){3}{\line(0,-1){1}}
\multiput(77.5,16)(-2,0){12}{\line(-1,0){1}}
\multiput(79,16.5)(0,2.4){3}{\line(0,1){1.2}}
\thinlines
\put(61,20){\line(-1,1){3}}
\put(67,18){\line(-1,1){6}}
\put(73,16){\line(-1,1){9}}
\put(79,14){\line(-1,1){12}}
\put(85,12){\line(-1,1){15}}
\put(90,11){\line(-1,1){17}}
\put(90,15){\line(-1,1){14}}
\put(88,21){\line(-1,1){9}}
\put(90,23){\line(-1,1){8}}
\put(90,27){\line(-1,1){5}}
\put(90,31){\line(-1,1){2}}
\put(51.8,25.7){$b_1$}
\put(51.3,42){${f}$}
\put(51.8,14.5){$b_2$}
\put(51.8,18.7){$O$}
\put(68.7,19.5){$a_1$}
\put(77.8,22.8){$a_2$}
\put(88.2,18.5){${e}$}
\put(60,-56){\ris\label{r99}}

\end{picture}}

\end{picture}
\end{wrapfigure}
\textit{2) If for each \fo{i\in\lbrace 1,2, \dots, n\rbrace} we have ${a_i-b_i\sqrt{p}>0}$, then a rectangle with side ratio $z$ can be tiled by rectangles with side ratios $x_1$,~\dots,~$x_n$ if and only if $${\displaystyle z\in\left\lbrace  e+f\sqrt{p}\,\vert\, e,f\in\q, e>0, \frac{|f|}{e}\leqslant\underset{1\leqslant i\leqslant n}{\max}\frac{|b_i|}{a_i}\right\rbrace=:M(x_1, \dots, x_n).}$$
3) If for each \fo{i\in\lbrace 1,2, \dots, n\rbrace} we have ${a_i-b_i\sqrt{p}<0}$, then a rectangle with side ratio $z$ can be tiled by rectangles with side ratios $x_1$,~\dots,~$x_n$ if and only if 
$${\displaystyle z\in\left\lbrace  e+f\sqrt{p}\,\vert\, e,f\in\q, f>0, \frac{|e|}{f}\leqslant\underset{1\leqslant i\leqslant n}{\max}\frac{|a_i|}{b_i}\right\rbrace=:N(x_1, \dots, x_n).}$$}

The assertions 2) and 3) of Theorem~1 have a simple geometric meaning. We represent the number of the form ${z=e+f\sqrt{p}}$ as the point~\fo{(e,f)} on the coordinate plane. In Figure~\ref{r99}.a the case when \fo{n=2}, ${a_1-b_1\sqrt{p}>0}$, and ${a_2-b_2\sqrt{p}>0}$ is shown: the shaded region depicts the set \fo{M(x_1,x_2)}. This is the least angle that is symmetrical with respect to the $e$-axis and contains all the points~\fo{(a_i,b_i)} (\fo{1\leqslant i\leqslant n}). Similarly, in Figure~\ref{r99}.b the case when \fo{n=2}, ${a_1-b_1\sqrt{p}<0}$, and ${a_2-b_2\sqrt{p}<0}$ is shown. The set \fo{N(x_1,x_2)} is shaded there. It is symmetric with respect to the $f$-axis.

It is clear that if for each \fo{i\in\lbrace 1,2, \dots, n\rbrace} we have ${a_i-b_i\sqrt{p}>0}$, then there exists a number \fo{k\in\lbrace 1,2, \dots, n\rbrace} such that \fo{M(x_1, \dots, x_n)=M(x_k)}. In this case, by the main theorem, if a rectangle can be tiled by rectangles with side ratios $x_1$,~\dots,~$x_n$, then it also can be tiled by rectangles with side ratio $x_k$. If for each \fo{i\in\lbrace 1,2, \dots, n\rbrace} we have ${a_i-b_i\sqrt{p}<0}$, then a similar assertion holds for \fo{N(x_1, \dots, x_n)}. 
\newpage
\glava\textbf{ Previous work.}\end{center}

Many mathematicians worked on the discussed problem (see [1--12]). Some of the obtained results are given below.

\teo\textbf{ (Dehn, 1903, \cite{D}).} \textit{A rectangle can be tiled by squares (not necessarily equal) if and only if the side ratio of the rectangle is rational.}\label{teorDm}

\teo\textbf{ (Dehn, 1903, \cite[Theorem 1]{FLR}).} \textit{If a rectangle with side ratio $x$ can be tiled by rectangles with side ratios \fo{x_1, x_2, \ldots, x_n}, then~$x$ can be written as a rational function with rational coefficients in the arguments \fo{x_1, x_2, \ldots, x_n}.}\label{teorD}


In spite of the simplicity of the statement of Theorem~\ref{teorDm}, the proof was initially quite complicated.

One of the most powerful methods of the proof of the impossibility of tilings is a physical interpretation using electrical networks. This method was discovered by R.~L.~Brooks, C.~A.~B.~Smith, A.~H.~Stone, and W.~T.~Tutte~\cite{BSST}. They successfully applied it to the proof of Theorem~\ref{teorDm} and for the search of dissections of squares into squares with different sides. Also the method was used, for example, in~\cite{PS}.

There are elementary proofs of Theorems~\ref{teorDm} and~\ref{teorD}: Theorem~\ref{teorDm} is proved in~\cite[Dehn's theorem]{KV1} and Theorem~\ref{teorD} is proved in~\cite{sobe_kva}.


\teo\textbf{ (Laczkovich, Rinne, Szekeres, Freiling, 1994, \cite{FR,LS}).} \textit{For a number \fo{r>0} the following 3 conditions are equivalent:\\
1) a square can be tiled by rectangles with side ratio $r$;\\
2) for certain positive rational numbers \fo{c_i} we have $$c_1r+\cfrac{1}{c_2r+\cfrac{1}{c_3r+\ldots+\cfrac{1}{c_nr}}}=1;$$
3) there exists a nonzero polynomial $P$ with integer coefficients such that the number $r$ is a root of $P$ and each complex root of $P$ has a positive real part.}\label{FRLS}

The following theorem is a particular case of Theorem~\ref{FRLS}, which has an elementary proof.

\teo\textbf{ (\cite{sobe_kva}).} \textit{Suppose \fo{a,b\in\q} and \fo{x=a+b\sqrt{2}>0}; then a square can be tiled by rectangles with side ratio $x$ if and only if \fo{a-b\sqrt{2}>0}.}



In~\cite{FLR} Freiling, Laczkovich, and Rinne reduced the problem of dissections of
a rectangle into rectangles similar to a given one to a hard algebraic
problem. They found an algebraic criterion of the existence of
a dissection. However, the criterion does not provide an algorithm for the construction of the required dissection. But they completely
solved the problem posed at the beginning of this paper in a particular case, when the side ratios of the rectangles are
quadratic irrationalities. We state their theorem (Theorem 7
in~\cite{FLR}) in the form equivalent to the original one.

\teo\textbf{ (Freiling, Laczkovich, Rinne, 1997, \cite[Theorem 7]{FLR}).} \textit{Suppose \fo{\alpha,\beta,\gamma,\delta,p\in\q}, \fo{\beta\sqrt{p}\notin\q}, ${u=\alpha+\beta\sqrt{p}>0}$, and ${v=\delta u+\gamma}$; then a rectangle with side ratio $v$ can be tiled by rectangles with side ratio $u$ if and only if either ${\gamma=0}$ and ${\delta>0}$ or else \fo{\alpha\neq 0}, \fo{\frac{\gamma(\alpha^2-\beta^2p)}{\alpha}>0}, and \fo{\delta+\frac{\gamma}{2\alpha}\geqslant 0}.}\label{TeoFLR}

This theorem solves the posed problem in the particular case when $n=1$.

Freiling, Laczkovich, and Rinne proved Theorem~\ref{TeoFLR} using the method of \2areas\3. In this paper we also use a similar method but in a simpler form.


In~\cite{KeKi} K.~Keating and J.~L.~King solved a related problem about \2signed\3 tilings. And in~\cite{Szegedy} B.~Szegedy solved a related problem on tilings of the square by similar right triangles.

\glava\textbf{ Definitions.}\end{center}


\opred Let $a$ and $b$ be positive real numbers; then each of the numbers $a/b$ and $b/a$ is called a \textit{side ratio} of a rectangle with sides $a$ and $b$.

\opred \fo{\R^+:=\left\lbrace r\in\R\,\vert\, r>0 \right\rbrace}, \fo{\R^-:=\left\lbrace r\in\R\,\vert\, r<0 \right\rbrace}, \fo{\q^+:=\left\lbrace q\in\q\,\vert\, q>0 \right\rbrace}, \fo{\q^+_0:=\q^+\cup\left\lbrace 0 \right\rbrace}, \fo{\q^-:=\left\lbrace q\in\q\,\vert\, q<0 \right\rbrace}.

\opred A real number $x$ is called a \textit{quadratic irrationality} if there exist rational numbers $a$, $b$, $p$ such that \fo{x=a+b\sqrt{p}} and \fo{\sqrt{p}\notin\q}. Denote the set of all quadratic irrationalities \fo{a+b\sqrt{p}} with fixed $p$ by \fo{\hor{p}}. And denote the set of all quadratic irrationalities \fo{a+b\sqrt{p}} with positive $a$ and $b$ and with fixed $p$ by \fo{\horp{p}}.

\opred $B(x_1, \dots, x_n):=\lbrace z\in\R^+\,\vert\, \mbox{a rectangle with side ratio }z\mbox{ can be tiled by rectangles with}$ $\mbox{side ratios }x_1, \dots, x_n\rbrace$.



\glava\textbf{ Lemmas.}\end{center}



\lemma \textit{If \fo{a,b\in B(x_1, \dots, x_n)}, then \fo{(a+b)\in B(x_1, \dots, x_n)}.}\label{l1}

$\blacktriangleleft\,$Suppose that \fo{a,b\in B(x_1, \dots, x_n)}; then a rectangle with sides 1~and~$a$ and a rectangle with sides 1~and~$b$ can be tiled by rectangles with side ratios~$x_1$,~\dots,~$x_n$. It is clear that then a rectangle with sides 1~and~$a+b$ can be tiled by rectangles with side ratios~$x_1$,~\dots,~$x_n$, i.e., \fo{(a+b)\in B(x_1, \dots, x_n)}.$\,\blacktriangleright$

\lemma\label{l2} \textit{If \fo{a\in B(x_1, \dots, x_n)}, then \fo{a^{-1}\in B(x_1, \dots, x_n)}.}

$\blacktriangleleft\,$By definition, if $a$ is a side ratio of a rectangle, then $a^{-1}$ is also a side ratio of the rectangle.$\,\blacktriangleright$

\lemma\label{l3} \textit{If \fo{a\in B(x_1, \dots, x_n)} and \fo{q\in\q^+}, then \fo{qa\in B(x_1, \dots, x_n)}.}

\begin{wrapfigure}[10]{l}{0.25\linewidth} 
\vspace{-6.5ex}
\begin{picture}(45,45) 
\put(5,33){\line(1,0){15}} 
\put(5,29){\line(1,0){15}}
\put(5,25){\line(1,0){15}}
\put(26,33){\line(1,0){15}} 
\put(26,29){\line(1,0){15}}
\put(26,25){\line(1,0){15}}
\put(5,5){\line(1,0){15}} 
\put(5,9){\line(1,0){15}}
\put(5,13){\line(1,0){15}}
\put(26,5){\line(1,0){15}} 
\put(26,9){\line(1,0){15}}
\put(26,13){\line(1,0){15}}
\put(5,33){\line(0,-1){10}}
\put(11,33){\line(0,-1){10}}
\put(17,33){\line(0,-1){10}}
\put(29,33){\line(0,-1){10}}
\put(35,33){\line(0,-1){10}}
\put(41,33){\line(0,-1){10}}
\put(5,15){\line(0,-1){10}}
\put(11,15){\line(0,-1){10}}
\put(17,15){\line(0,-1){10}}
\put(29,15){\line(0,-1){10}}
\put(35,15){\line(0,-1){10}}
\put(41,15){\line(0,-1){10}}
\put(21,33){\circle*{0.5}}
\put(23,33){\circle*{0.5}}
\put(25,33){\circle*{0.5}}
\put(21,5){\circle*{0.5}}
\put(23,5){\circle*{0.5}}
\put(25,5){\circle*{0.5}}
\put(5,17){\circle*{0.5}}
\put(5,19){\circle*{0.5}}
\put(5,21){\circle*{0.5}}
\put(41,17){\circle*{0.5}}
\put(41,19){\circle*{0.5}}
\put(41,21){\circle*{0.5}}
\put(3.5,4.5){$^1$}
\put(3.5,8.5){$^1$}
\put(3.5,24.5){$^1$}
\put(3.5,28.5){$^1$}
\put(7.5,32){$^a$}
\put(13.5,32){$^a$}
\put(37.5,32){$^a$}
\put(31.5,32){$^a$}
\put(20,2.1){$ma$}
\put(41.5,18.5){$n$}
\put(17,-4){\ris\label{ri2}}
\end{picture}
\end{wrapfigure}

$\blacktriangleleft\,$Suppose that~\fo{a\in B(x_1, \dots, x_n)} and~${q=\frac{m}{n}}$, where~\fo{m,n\in\n}. A rectangle with sides~$n$ and~$ma$ can be tiled by rectangles with sides~1 and~$a$; see Figure~\ref{ri2}. Also, the rectangles with sides~1 and~$a$ can be tiled by rectangles $x_1$,~\dots,~$x_n$, because \fo{a\in B(x_1, \dots, x_n)}. Therefore, ${\frac{ma}{n}=qa\in B(x_1, \dots, x_n)}$.$\,\blacktriangleright$
 
\lemma\label{l4} \textit{Suppose \fo{x:=a+b\sqrt{p}}, \fo{\overline{x}:=a-b\sqrt{p}}, and \fo{x,\x\in\hor{p}\cap\R^+}; then \fo{\x\in B(x)} and \fo{\q^+\subset B(x)}.}

$\blacktriangleleft\,$First we prove that \fo{\x\in B(x)}. Using Lemma~\ref{l2}, we get $${\frac{a-b\sqrt{p}}{a^2-pb^2}=\frac{1}{a+b\sqrt{p}}=\frac{1}{x}\in B(x).}$$ By Lemma~\ref{l3}, it follows that $${\x=a-b\sqrt{p}=\frac{a-b\sqrt{p}}{a^2-pb^2}(a^2-pb^2)\in B(x)}$$ because \fo{(a^2-pb^2)\in\q^+}. 

Now we prove that \fo{\q^+\subset B(x)}. Using Lemma~\ref{l1}, we get \fo{2a=x+\x\in B(x)} because \fo{x,\x\in B(x)}. By Lemma~\ref{l3}, it follows that for any \fo{q\in\q^+} we have \fo{2aq\in B(x)}, i.e., \fo{\q^+\subset B(x)}.$\,\blacktriangleright$

\lemma\label{l5} \textit{Suppose \fo{x:=a+b\sqrt{p}\in\hor{p}\cap\R^+} and \fo{\overline{x}:=a-b\sqrt{p}\in\hor{p}\cap\R^-}; then \fo{(-\x)\in B(x)} and for each \fo{q\in\q^+} we have \fo{q\sqrt{p}\in B(x)}.}

The proof is analogous to the proof of Lemma~\ref{l4} and is left to the reader.

\lemma\label{l7} \textit{If~\fo{x_1, \dots, x_n\in\hor{p}\cap\R^+}, then~\fo{B(x_1, \dots, x_n)\subset\hor{p}\cap\R^+}.} 

$\blacktriangleleft\,$Suppose~\fo{x_1, \dots, x_n\in\hor{p}\cap\R^+} and \fo{z\in B(x_1, \dots, x_n)}. Then by Theorem~\ref{teorD} the number $z$ can be written as a rational function with rational coefficients of the arguments \fo{x_1, x_2, \ldots, x_n}. Hence \fo{z\in\hor{p}}. Furthermore, by definition, \fo{z>0}. Therefore, \fo{B(x_1, \dots, x_n)\subset\hor{p}\cap\R^+}.$\,\blacktriangleright$

If \fo{z\in M(x_1, \dots, x_n)} or \fo{z\in N(x_1, \dots, x_n)}, then \fo{z\in\hor{p}\cap\R^+}. Also, if \fo{z\in B(x_1, \dots, x_n)}, then \fo{z\in\hor{p}\cap\R^+}. Therefore it suffices to prove Theorem 1 in the case when the condition \fo{z\in\hor{p}\cap\R^+} holds. In what follows we assume that \fo{e+f\sqrt{p}=z\in\hor{p}\cap\R^+}.

\opred Let \fo{A,B,C\in\R}. The \textit{``area''} of a rectangle with sides ${\alpha+\beta\sqrt{p}}$ and ${\gamma+\delta\sqrt{p}}$, where \fo{\alpha,\beta,\gamma,\delta,p\in\q} and \fo{\sqrt{p}\notin\q}, is the number $$S:=\alpha\gamma A+\beta\gamma B+\alpha\delta B+\beta\delta C.$$\label{o5}

It is easy to prove that the definition of the \2area\3 of a rectangle is invariant under the permutation of two sides of the rectangle.

\lemma\label{l9} \textit{If a rectangle $R$ is tiled by rectangles \fo{R_1, \dots, R_n} with the sides from \fo{\hor{p}}, then the \2area\3 of $R$ is equal to the sum of the \2areas\3 of \fo{R_1, \dots, R_n}.}

\begin{wrapfigure}[7]{l}{0.25\linewidth} 
\vspace{0ex}
\begin{picture}(45,83)
\thicklines
\put(0,80){\line(1,0){45}}
\put(14,72){\line(1,0){31}}
\put(21,76){\line(1,0){24}}
\put(0,60){\line(1,0){45}}
\put(0,60){\line(0,1){20}}
\put(14,60){\line(0,1){20}}
\put(21,72){\line(0,1){8}}
\put(25,60){\line(0,1){12}}
\put(35,60){\line(0,1){12}}
\put(45,60){\line(0,1){20}}
\thinlines
\multiput(0,72)(2,0){8}{\line(1,0){1}}
\multiput(0,76)(2,0){11}{\line(1,0){1}}
\multiput(21,60)(0,2){6}{\line(0,1){1}}
\multiput(25,80)(0,-2){5}{\line(0,-1){1}}
\multiput(35,80)(0,-2){5}{\line(0,-1){1}}
\put(15.5,55){\ris\label{rrr5}}
\end{picture}
\end{wrapfigure}
$\blacktriangleleft\,$Let $P_1$ and $P_2$ be two rectangles with the sides from \fo{\hor{p}} having a common side. Then it is clear that the sum of the \2areas\3 of $P_1$ and $P_2$ is equal to the \2area\3 of the rectangle \fo{P_1\cup P_2}.

For arbitrary $n$, we extend each cut as much as possible; see Figure~\ref{rrr5}. Then it is obvious that each rectangle of the resulting tiling has side ratio that is a quadratic irrationality. Using the assertion of the previous paragraph, we get the assertion of Lemma~\ref{l9}.$\,\blacktriangleright$

\glava\textbf{ Proof of the existence.}\end{center} 

In this section we prove the \2if\3 part in each of the assertions 1)--3) of Theorem~\ref{mainT}.

\textbf{Assertion 1).} Suppose that there exist \fo{i,j\in\lbrace 1,2, \dots, n\rbrace} such that \fo{\x_i:=a_i-b_i\sqrt{p}>0} and \fo{\x_j:=a_j-b_j\sqrt{p}<0}. We shall see that \fo{P:=\hor{p}\cap\R^+\subset B(x_1, \dots, x_n)}. Suppose ${z=e+f\sqrt{p}\in P}$. Let us consider 3 cases.



Case 1:~\fo{e>0},~\fo{f>0}. We have \fo{x_i,\x_i\in\hor{p}\cap\R^+} and \fo{e\in\q^+}. Then, using Lemma~\ref{l4}, we get \fo{e\in B(a_i+b_i\sqrt{p})=B(x_i)\subset B(x_1, \dots, x_n)}. We also have \fo{x_j\in\hor{p}\cap\R^+}, \fo{\x_j\in\hor{p}\cap\R^-}, and \fo{f\in\q^+}. Using Lemma~\ref{l5}, we get \fo{f\sqrt{p}\in B(a_j+b_j\sqrt{p})=B(x_j)\subset B(x_1, \dots, x_n)}. Hence, by Lemma~\ref{l1}, we have \fo{e+f\sqrt{p}\in B(x_1, \dots, x_n)}. 

Case 2:~\fo{e<0} (obviously, in this case \fo{f>0} and \fo{pf^2-e^2>0}). Then ${\frac{1}{z}=\frac{-e+f\sqrt{p}}{pf^2-e^2}\in\horp{p}\subset B(x_1, \dots, x_n)}$ by Case~1. Using Lemma~\ref{l2}, we get \fo{z\in B(x_1, \dots, x_n)}.

Case 3:~\fo{f<0} (obviously, in this case \fo{e>0} and \fo{e^2-pf^2>0}). Then ${\frac{1}{z}=\frac{e-f\sqrt{p}}{e^2-pf^2}\in\horp{p}\subset B(x_1, \dots, x_n)}$ by Case~1. Using Lemma~\ref{l2}, we get \fo{z\in B(x_1, \dots, x_n)}.



\textbf{Assertion 2).} Assume without loss of generality that ${|b_1|/a_1=\underset{1\leqslant i\leqslant n}{\max}(|b_i|/a_i)}$ for each \fo{i\in\lbrace 1, \dots, n\rbrace}. If~\fo{b_1=0}, then all the numbers $x_1$,~\dots,~$x_n$ are rational, hence \fo{M(x_1, \dots, x_n)=\q^+}. Using Lemma~\ref{l3}, we get \fo{M(x_1, \dots, x_n)=\q^+\subset B(x_1, \dots, x_n)} . In what follows we assume that \fo{b_1\neq 0}.


Suppose~\fo{z=(e+f\sqrt{p})\in M(x_1, \dots, x_n)}. By the conditions of assertion~2) of Theorem~\ref{mainT}, we have ${e\in\q^+}$, ${f\in\q}$, and ${|f|/e\leqslant |b_1|/a_1}$. Let us consider two cases.

Case 1: ${f\in\q^+_0}$. Let us prove that $$z=e+f\sqrt{p}=\frac{f}{|b_1|}(a_1+|b_1|\sqrt{p})+\left( e-\frac{fa_1}{|b_1|}\right)\in B(x_1, \dots, x_n).$$

Indeed, ${\left( e-\frac{fa_1}{|b_1|}\right)\in\q^+_0}$ (because~${\frac{|f|}{e}\leqslant\frac{|b_1|}{a_1}}$) and~${\frac{f}{|b_1|}\in\q^+_0}$ (because~${f\in\q^+_0}$). It follows from Lemma~\ref{l4} that ${(a_1-b_1\sqrt{p})\in B(x_1, \dots, x_n)}$, hence ${(a_1+|b_1|\sqrt{p})\in B(x_1, \dots, x_n)}$. It also follows from Lemma~\ref{l4} that for each \fo{q\in\q^+} we have \fo{q\in B(x_1, \dots, x_n)}. Using Lemma~\ref{l1} and Lemma~\ref{l3}, we get the required inclusion.

Case 2:~${f\in\q^-}$. Let us prove that $$z=e+f\sqrt{p}=\frac{|f|}{|b_1|}(a_1-|b_1|\sqrt{p})+\left( e-\frac{|f|a_1}{|b_1|}\right)\in B(x_1, \dots, x_n).$$

Indeed,~${\left( e-\frac{|f|a_1}{|b_1|}\right)\in\q^+_0}$ (because ${\frac{|f|}{e}\leqslant\frac{|b_1|}{a_1}}$) and ${\frac{|f|}{|b_1|}\in\q^+}$. Using Lemma~\ref{l1} and Lemma~\ref{l3}, we get the required inclusion.

The proof of the \2if\3 part of Assertion 3) of Theorem~\ref{mainT} is analogous to the proof of the \2if\3 part of Assertion~2).


\glava\textbf{ Proof of the impossibility.}\end{center}

In this section we prove the \2only if\3 part in each of the assertions 1)--3) of Theorem~\ref{mainT}.

\textbf{Assertion 1).} By Lemma~\ref{l7}, it follows that~\fo{B(x_1, \dots, x_n)\subset\hor{p}\cap\R^+}.

\textbf{Assertions 2) and 3).} First we prove the \2only if\3 part in the particular case when \fo{n=1}. In other words, we prove that if \fo{a_1-b_1\sqrt{p}>0}, then \fo{z\notin M(x_1)\Rightarrow z\notin B(x_1)}; and if \fo{a_1-b_1\sqrt{p}<0}, then \fo{z\notin N(x_1)\Rightarrow z\notin B(x_1)}.

Assume the converse, i.e., assume that \fo{z\in B(x_1)} and one of the following conditions holds:\\•~\fo{a_1-b_1\sqrt{p}>0} and \fo{|f|/e>|b_1|/a_1};\\•~\fo{a_1-b_1\sqrt{p}<0} and \fo{|e|/f>|a_1|/b_1}.

The following two definitions are useful for sequel. The rectangle with side ratio~$z$ is called \textit{great} and the rectangles with side ratio~$x_1$ are called~\textit{small}.

Let the great rectangle have sides 1 and \fo{e+f\sqrt{p}}. By the theorem on electric current and lengths of sides and by the lemma on the resistance of an electrical circuit from~\cite{sobe_kva}, the sides of the small rectangles can be expressed through the number $x_1$ using addition, subtraction, multiplication, and division. Then, since the set \fo{\hor{p}} is closed with respect to these operations, the lengths of the sides of the small rectangles are quadratic irrationalities. Thus, we can use Definition~5 for the \2areas\3 of the great and small rectangles.
Set $$A=f,\,\,\,B=-e,\,\,\,C=\frac{2fa_1^2}{b_1^2}-pf.$$ Then the \2area\3 of the great rectangle in this case is $$S_{\mbox{\textit{\scriptsize{G}}}}=e\cdot1\cdot f+f\cdot1\cdot(-e)+e\cdot0\cdot(-e)+(-e)\cdot0\cdot\left(\frac{2fa_1^2}{b_1^2}-pf\right)=0.$$

Now we estimate the \2area\3 of the small rectangles. Let the sides of some of the small rectangles is ${\alpha+\beta\sqrt{p}}$ and ${\gamma+\delta\sqrt{p}}$. Then \fo{\frac{\gamma+\delta\sqrt{p}}{\alpha+\beta\sqrt{p}}=a_1+b_1\sqrt{p}}. Then \fo{\gamma+\delta\sqrt{p}=\alpha a_1+p\beta b_1+\sqrt{p}(\beta a_1+\alpha b_1)} and the following equalities hold: $$\gamma=\alpha a_1+p\beta b_1\,\,\,\mbox{ and }\,\,\,\delta=\beta a_1+\alpha b_1.$$ Therefore the \2area\3 of the small rectangle is written as $$S_{\mbox{\textit{\scriptsize{S}}}}=\alpha\gamma f-\beta\gamma e-\alpha\delta e+\beta\delta\left(\frac{2fa_1^2}{b_1^2}-pf\right)=$$
$$=\alpha f(\alpha a_1+p\beta b_1)-\beta e(\alpha a_1+p\beta b_1)-\alpha e(\beta a_1+\alpha b_1)+\beta\left(\frac{2fa_1^2}{b_1^2}-pf\right)(\beta a_1+\alpha b_1)=$$
$$=\alpha^2(fa_1-eb_1)+2\alpha\beta(\frac{fa_1^2}{b_1}-ea_1)+\beta^2(\frac{2fa_1^3}{b_1^2}-pfa_1-peb_1).$$

Note that by the assumption we have \fo{fa_1-eb_1\neq 0}. So the expression \fo{S_{\mbox{\textit{\scriptsize{S}}}}/\beta^2} is a square trinomial in \fo{\alpha/\beta}. Let us show that its discriminant $D$ is negative. This will implies that for fixed numbers $a_1$, $b_1$, $e$, $f$ the value of \fo{S_{\mbox{\textit{\scriptsize{S}}}}} either always positive or always negative, hence the sum of all the \2areas\3 of small rectangles is not equal to zero, i.e., not equal to the \2area\3 of the great rectangle. Thus, we shall get a contradiction.

So $$\frac{D}{4}=\frac{f^2a_1^4}{b_1^2}-\frac{2fea_1^3}{b_1}+e^2a_1^2-\frac{2f^2a_1^4}{b_1^2}+pf^2a_1^2+pfa_1eb_1+\frac{2fea_1^3}{b_1}-pfa_1eb_1-pe^2b_1^2=$$
$$=\frac{-f^2a_1^4}{b_1^2}+e^2a_1^2+pf^2a_1^2-pe^2b_1^2=(a_1^2-pb_1^2)\left( e^2-\frac{f^2a_1^2}{b_1^2}\right).$$ We have two cases:\\•~\fo{a_1-b_1\sqrt{p}<0} and \fo{|e|/f>|a_1|/b_1\Leftrightarrow |e|>f|a_1|/b_1};\\•~\fo{a_1-b_1\sqrt{p}>0} and \fo{|f|/e>|b_1|/a_1\Leftrightarrow e<a_1|f|/|b_1|}.\\In both cases \fo{D<0}.

The \2only if\3 part in the particular case when \fo{n=1} is proved. Above we proved the \2if\3 part of Theorem~\ref{mainT}. Therefore we have already proved the whole Theorem~\ref{mainT} in the particular case when \fo{n=1}.

Now let us prove the \2only if\3 part in Assertion~2) for arbitrary $n$. Assume without loss of generality that ${|b_1|/a_1=\underset{1\leqslant i\leqslant n}{\max}(|b_i|/a_i)}$. Let us prove that then \fo{B(x_1)=B(x_1, \dots, x_n)}.

We first prove that \fo{B(x_1, \dots, x_n)\subset B(x_1)}. Take \fo{z\in B(x_1, \dots, x_n)}. In a dissection of a rectangle with side ratio~$z$ into rectangles with side ratios of~$x_1$,~\dots,~$, x_n$, replace each rectangle with side ratio belonging to the set \fo{\lbrace x_2, \dots, x_n\rbrace} with its dissection into rectangles with side ratio~$x_1$ (we can do this, because by Theorem~\ref{mainT} in the case when \fo{n=1}, we have \fo{x_2, \dots, x_n\in B(x_1)}). We get a dissection of the rectangle with side ratio~$z$ into rectangles with side ratio~$x_1$. Thus \fo{B(x_1, \dots, x_n)\subset B(x_1)}.

The reverse inclusion is obvious from the definition of the set \fo{B(x_1, \dots, x_n)}. Thus \fo{B(x_1)=B(x_1, \dots, x_n)}.

By Theorem~\ref{mainT} in the case when \fo{n=1}, we get $${\displaystyle B(x_1, \dots, x_n)=B(x_1)=\left\lbrace e+f\sqrt{p}\,\vert\, e\in\q^+, f\in\q,\frac{|f|}{e}\leqslant\frac{|b_1|}{a_1} \right\rbrace.}$$ By the assumption that ${|b_1|/a_1=\underset{1\leqslant i\leqslant n}{\max}(|b_i|/a_i)}$ we get $${\displaystyle B(x_1, \dots, x_n)=\left\lbrace e+f\sqrt{p}\,\vert\, e\in\q^+, f\in\q,\frac{|f|}{e}\leqslant\underset{1\leqslant i\leqslant n}{\max}\frac{|b_i|}{a_i} \right\rbrace,}$$ that is required.

The proof of the \2only if\3 part of Assertion 3) of Theorem~\ref{mainT} is analogous to the proof of the \2only if\3 part of Assertion~2).

\textbf{Acknowledgements.} The author is grateful to assistant professor of HSE Mikhail Skopenkov for constant attention to this work.

\bibliographystyle{elsarticle-num}

\end{document}


\begin{center}
{\Large \textbf{Разрезания прямоугольника на прямоугольники с заданными отношениями сторон}}\\
$\phantom{kl}$\\
\begin{large}
\textbf{Фёдор Шаров\footnote{Автор был частично поддержан грантом президента Российской Федерации MK-6137.2016.1.}}\\
$\phantom{kl}$\\
\end{large}
{\small\textit{Факультет математики, Национальный исследовательский университет Высшая школа экономики}}
\end{center}

\glava\textbf{ Основной результат.}\end{center}

В этой статье мы рассмотрим такую задачу: \textit{дан набор прямоугольников; какие фигуры можно разрезать на прямоугольники, подобные данным?} Поставленная задача в целом сложна, до сих пор не имеет решения и вряд ли может быть в разумном смысле решена в общем виде. Однако, довольно интересным является исследование частных случаев этой общей проблемы. В нашей статье мы выясним, какие прямоугольники можно разрезать на подобные $n$ данным при условии, что отношения сторон данных прямоугольников --- квадратичные иррациональности.


Сформулируем основную теорему нашей статьи. Эта теорема является новым результатом, который ранее не был опубликован. При этом доказательство использует только элементарные методы. В <<Математическом просвещении>>, как правило, не публикуются исследовательские статьи, но в тех случаях, когда речь идёт об общедоступных статьях, бывают исключения (см., например,~\cite{BiTk}). По нашему мнению, эта статья как раз относится к такой категории.

\teo\textbf{ (основная).} \textit{Пусть ${x_1=a_1+b_1\sqrt{p}}$,~\dots,~${x_n=a_n+b_n\sqrt{p}}$ --- такие числа, что ${x_i>0}$, $a_i,b_i,p\in\q$ (\fo{1\leqslant i\leqslant n}) и ${\sqrt{p}\notin\q}$. Тогда:\\
1) если существуют такие числа $i$ и $j$, что \fo{1\leqslant i,j\leqslant n} и ${(a_i-b_i\sqrt{p})(a_j-b_j\sqrt{p})<0}$, то прямоугольник с отношением сторон $z$ можно разрезать на прямоугольники с отношениями сторон $x_1$,~\dots,~$x_n$ тогда и только тогда, когда}

\begin{wrapfigure}[22]{l}{0.27\linewidth} 
\vspace{-5.3ex}
\hspace*{2.75ex}
\begin{picture}(45,100) 
\thicklines
\put(-5,86){а)}
\put(-5,41){б)}
\put(0,55){\vector(1,0){42}}
\put(22,52){\vector(0,1){38}}
\put(22,55){\vector(1,3){5}}
\put(22,55){\vector(-1,4){6}}
\multiput(27.3,70.9)(2.4,7.2){3}{\line(1,3){1.2}}
\multiput(21.7,55.9)(-2.4,7.2){5}{\line(-1,3){1.2}}
\multiput(27,69.5)(0,-2){8}{\line(0,-1){1}}
\multiput(26.5,70)(-2,0){3}{\line(-1,0){1}}
\multiput(16,77.5)(0,-2){12}{\line(0,-1){1}}
\multiput(16.5,79)(2.4,0){3}{\line(1,0){1.2}}
\thinlines
\put(20,61){\line(1,-1){3}}
\put(18,67){\line(1,-1){6}}
\put(16,73){\line(1,-1){9}}
\put(14,79){\line(1,-1){12}}
\put(12,85){\line(1,-1){15}}
\put(11,90){\line(1,-1){17}}
\put(15,90){\line(1,-1){14}}
\put(21,88){\line(1,-1){9}}
\put(23,90){\line(1,-1){8}}
\put(27,90){\line(1,-1){5}}
\put(31,90){\line(1,-1){2}}
\put(25.8,52){$a_1$}
\put(40.4,51.4){$e$}
\put(14.8,52){$a_2$}
\put(18.8,51.7){$O$}
\put(18.9,68.3){$b_1$}
\put(22.5,77.9){$b_2$}
\put(18.4,87.6){$f$}

\put(-48,0){\begin{picture}(100,45)

\thicklines
\put(55,0){\vector(0,1){45}}
\put(52,22){\vector(1,0){38}}
\put(55,22){\vector(3,1){15}}
\put(55,22){\vector(4,-1){24}}
\multiput(70.9,27.3)(7.2,2.4){3}{\line(3,1){4}}
\multiput(55.9,21.7)(7.2,-2.4){5}{\line(3,-1){4}}
\multiput(69.5,27)(-2,0){8}{\line(-1,0){1}}
\multiput(70,26.5)(0,-2){3}{\line(0,-1){1}}
\multiput(77.5,16)(-2,0){12}{\line(-1,0){1}}
\multiput(79,16.5)(0,2.4){3}{\line(0,1){1.2}}
\thinlines
\put(61,20){\line(-1,1){3}}
\put(67,18){\line(-1,1){6}}
\put(73,16){\line(-1,1){9}}
\put(79,14){\line(-1,1){12}}
\put(85,12){\line(-1,1){15}}
\put(90,11){\line(-1,1){17}}
\put(90,15){\line(-1,1){14}}
\put(88,21){\line(-1,1){9}}
\put(90,23){\line(-1,1){8}}
\put(90,27){\line(-1,1){5}}
\put(90,31){\line(-1,1){2}}
\put(51.8,25.7){$b_1$}
\put(51.3,42){${f}$}
\put(51.8,14.5){$b_2$}
\put(51.8,18.7){$O$}
\put(68.7,19.5){$a_1$}
\put(77.8,22.8){$a_2$}
\put(88.2,18.5){${e}$}
\put(60,-6){\ris\label{r99}}

\end{picture}}

\end{picture}
\end{wrapfigure}
$${\displaystyle z\in\left\lbrace  e+f\sqrt{p}>0\,\vert\, e,f\in\q \right\rbrace;}$$
\textit{2) если для всех $i$ (\fo{1\leqslant i\leqslant n}) ${a_i-b_i\sqrt{p}>0}$, то прямоугольник с отношением сторон $z$ можно разрезать на прямоугольники с отношениями сторон $x_1$,~\dots,~$x_n$ тогда и только тогда, когда $${\displaystyle z\in\left\lbrace  e+f\sqrt{p}\,\vert\, e,f\in\q, e>0, \frac{|f|}{e}\leqslant\underset{i}{\max}\frac{|b_i|}{a_i}\right\rbrace=:M(x_1, \dots, x_n);}$$
3) если для всех $i$ (\fo{1\leqslant i\leqslant n}) ${a_i-b_i\sqrt{p}<0}$, то прямоугольник с отношением сторон $z$ можно разрезать на прямоугольники с отношениями сторон $x_1$,~\dots,~$x_n$ тогда и только тогда, когда 
$${\displaystyle z\in\left\lbrace  e+f\sqrt{p}\,\vert\, e,f\in\q, f>0, \frac{|e|}{f}\leqslant\underset{i}{\max}\frac{|a_i|}{b_i}\right\rbrace=:N(x_1, \dots, x_n).}$$}

Условия второго и третьего пунктов теоремы имеют простой геометрический смысл. Число вида ${z=e+f\sqrt{p}}$ изобразим на координатной плоскости в виде точки~\fo{(e,f)}. На рис.~\ref{r99}.а показан случай, когда \fo{n=2}, ${a_1-b_1\sqrt{p}<0}$ и ${a_2-b_2\sqrt{p}<0}$: заштрихованная область изображает множество \fo{N(x_1,x_2)}. Это наименьший симметричный относительно оси~\fo{Of} угол, который содержит все точки~\fo{(a_i,b_i)} (\fo{1\leqslant i\leqslant n}). Аналогично на рис.~\ref{r99}.б показан случай, когда \fo{n=2}, ${a_1-b_1\sqrt{p}>0}$ и ${a_2-b_2\sqrt{p}>0}$: здесь изображено множество \fo{M(x_1,x_2)}, симметричное относительно оси \fo{Oe}.

Ясно, что если ${a_i-b_i\sqrt{p}>0}$ (\fo{1\leqslant i\leqslant n}), то множество \fo{M(x_1, \dots, x_n)} всегда совпадает с одним из множеств $M(x_1)$,~\dots,~$M(x_n)$ --- пусть \fo{M(x_1, \dots, x_n)=M(x_k)} (\fo{k\in\lbrace 1, \dots, n\rbrace}). В этом случае, согласно основной теореме, все прямоугольники, которые можно разрезать на прямоугольники с отношениями сторон $x_1$,~\dots,~$x_n$, можно разрезать и на прямоугольники с отношением сторон $x_k$. Аналогичное утверждение справедливо и для множества \fo{N(x_1, \dots, x_n)}, когда ${a_i-b_i\sqrt{p}<0}$ (\fo{1\leqslant i\leqslant n}).

\glava\textbf{ Обзор известных результатов.}\end{center}

Многие математики занимались обсуждаемой проблемой (см. [1--14]). Ниже приведены несколько из полученных результатов.

\teo\textbf{ (Ден, 1903, см. \cite{D}).} \textit{Если прямоугольник разрезан на квадраты (не обязательно равные), то отношение его сторон рационально.}\label{teorDm}

\teo\textbf{ (Ден, 1903, см. теорему 1 в \cite{FLR}).} \textit{Если прямоугольник с отношением сторон~$x$ можно разрезать на прямоугольники с отношениями сторон \fo{x_1, x_2, \ldots, x_n}, то число~$x$ можно выразить через числа \fo{x_1, x_2, \ldots, x_n} с помощью сложения, вычитания, умножения и деления.}\label{teorD}

У теорем \ref{teorDm} и \ref{teorD} существуют элементарные доказательства: теорема \ref{teorDm} доказана в \cite[теорема Дена]{KV1}, теорема \ref{teorD} доказана в \cite{sobe_kva}, но не сформулирована там явно: её доказательство сразу следует из двух результатов: теоремы о сопротивлении цепи и отношении сторон прямоугольника и леммы о сопротивлении цепи.

\teo\textbf{ (Ласкович, Ринн, Секереш, Фрайлинг, 1994, см. \cite{FR,LS}).} \textit{Для числа \fo{r>0} следующие три условия эквивалентны:\\
1) квадрат можно разрезать на прямоугольники с отношением сторон $r$;\\
2) для некоторых положительных рациональных чисел \fo{c_i} выполнено равенство $$c_1r+\cfrac{1}{c_2r+\cfrac{1}{c_3r+\ldots+\cfrac{1}{c_nr}}}=1;$$
3) число $r$ является корнем ненулевого многочлена с целыми коэффициентами, у которого все комплексные корни имеют положительную действительную часть.}

Частным случаем теоремы Ласковича-Ринна-Секереша-Фрайлинга является следующая теорема, элементарное доказательство которой приведено в~\cite{sobe_kva}.

\teo\textbf{.} \textit{Пусть \fo{x=a+b\sqrt{2}>0}, где \fo{a,b\in\q}. Тогда из прямоугольников с отношением сторон $x$ можно составить квадрат тогда и только тогда, когда \fo{a-b\sqrt{2}>0}.}

В статье~\cite{FLR} Фрайлинг, Ласкович и Ринн свели задачу о разрезании прямоугольника на прямоугольники, подобные данному, к сложной алгебраической проблеме --- они дали алгебраический критерий возможности разрезания, правда, не дающий алгоритма проверки существования разрезания. Но зато это позволило им решить задачу для частного случая, когда отношения сторон являются квадратичными иррациональностями. Вот их теорема (теорема~7 в~\cite{FLR}) в формулировке, равносильной авторской.

\teo\textbf{ (Фрайлинг, Ласкович, Ринн, 1997).} \textit{Пусть ${u=\alpha+\beta\sqrt{p}>0}$, где \fo{\alpha,\beta,p\in\q} и \fo{\beta\sqrt{p}\notin\q}. И пусть ${v=\delta u+\gamma}$ для некоторых рациональных $\gamma$ и $\delta$. Тогда прямоугольник с отношением сторон $v$ можно разрезать на прямоугольники с отношением сторон $u$ тогда и только тогда, когда выполнено одно из двух условий:\\
1) ${\gamma=0}$ и ${\delta>0}$;\\ 
2) \fo{\alpha\neq 0}, \fo{\frac{\gamma(\alpha^2-\beta^2p)}{\alpha}>0} и \fo{\delta+\frac{\gamma}{2\alpha}\geqslant 0}.}\label{TeoFLR}

Эта теорема равносильна частному случаю основной теоремы при \fo{n=1} (равносильность доказана в главе 6). Поэтому в главе 5 попутно с доказательством основной теоремы мы получим и элементарное доказательство теоремы~\ref{TeoFLR}. 

К.~Китинг и Дж.Л.~Кинг решили близкую к поставленной в начале статьи задачу --- о разрезании прямоугольника на прямоугольники и так называемые <<антипрямоугольники>> (см.~\cite{KeKi,KK}).

\glava\textbf{ Определения и обозначения.}\end{center}

Для доказательства теоремы сначала дадим несколько определений и введём обозначения.

\opred Будем говорить, что прямоугольник $P$ \textit{разрезан} на прямоугольники $P_1$,~\dots,~$P_n$, если \fo{\overset{n}{\underset{i=1}{\cup}}P_i=P} и внутренности любых двух прямоугольников из $P_1$,~\dots,~$P_n$ не пересекаются. \textit{Разрезанием} прямоугольника $P$ на прямоугольники $P_1$,~\dots,~$P_n$ называется множество прямоугольников $P_1$,~\dots,~$P_n$ вместе с их фиксированным положением на плоскости относительно друг друга, если выполнено условие \fo{\overset{n}{\underset{i=1}{\cup}}P_i=P} и внутренности прямоугольников $P_1$,~\dots,~$P_n$ не пересекаются.

\opred \textit{Отношением сторон} прямоугольника со сторонами $a$ и $b$ будем называть каждое из чисел $a/b$ и $b/a$.

Множество всех положительных рациональных чисел будем обозначать $\q^+$, множество \fo{\q^+\cup\left\lbrace 0 \right\rbrace} будем обозначать \fo{\q^+_0}, множество всех отрицательных рациональных чисел будем обозначать $\q^-$. Аналогично для удобства введём обозначения \fo{\R^+:=\left\lbrace r\in\R\,\vert\, r>0 \right\rbrace} и \fo{\R^-:=\left\lbrace r\in\R\,\vert\, r<0 \right\rbrace}.

\opred Все числа, которые можно представить в виде \fo{x=a+b\sqrt{p}} с рациональным положительным $p$, таким, что $p$ не является квадратом рационального числа, и рациональными $a$ и $b$, называются \textit{квадратичными иррациональными числами (КИЧ)} или \textit{квадратичными иррациональностями}. Множество всех квадратичных иррациональностей \fo{x=a+b\sqrt{p}} при фиксированном $p$ будем обозначать $\hor{p}$. Множество всех КИЧ \fo{x=a+b\sqrt{p}} с положительными $a$ и $b$ при фиксированном $p$ будем обозначать $\horp{p}$.

\begin{wrapfigure}[11]{l}{0.29\linewidth} 
\vspace{-1.25ex}
\hspace*{2.75ex}
\begin{picture}(45,80) 
\put(-5,77){а)}
\put(0,80){\line(1,0){45}}
\put(15,72){\line(1,0){30}}
\put(21,76){\line(1,0){24}}
\put(0,60){\line(1,0){45}}
\put(0,60){\line(0,1){20}}
\put(15,60){\line(0,1){20}}
\put(21,72){\line(0,1){8}}
\put(25,60){\line(0,1){12}}
\put(35,60){\line(0,1){12}}
\put(45,60){\line(0,1){20}}
\put(-5,52){б)}
\put(0,55){\line(1,0){45}}
\put(0,35){\line(1,0){45}}
\put(0,48){\line(1,0){33}}
\put(15,42){\line(1,0){30}}
\put(0,35){\line(0,1){20}}
\put(15,35){\line(0,1){13}}
\put(33,42){\line(0,1){13}}
\put(45,35){\line(0,1){20}}
\put(16,30){\ris\label{rrr}}
\end{picture}
\end{wrapfigure}

\opred Любое разрезание прямоугольника на один прямоугольник будем называть \textit{тривиальным}. Далее тривиальные разрезания определяются по индукции: если даны два тривиальных разрезания двух прямоугольников \fo{P_1} и \fo{P_2} с общей стороной, но без общих внутренних точек, то разрезание их объединения, полученное объединением данных разрезаний прямоугольников \fo{P_1} и \fo{P_2}, также будем называть \textit{тривиальным}.

На рис. \ref{rrr}.а изображено тривиальное разрезание, на рис. \ref{rrr}.б --- нетривиальное. Неформально говоря, нужно представить разрезаемый прямоугольник в виде прямоугольного листа бумаги с нарисованным на нём разбиением на прямоугольники. Этот бумажный прямоугольник разрешается разрезать вдоль любого отрезка на два прямоугольника, потом производить такие операции по-отдельности с каждой из получившихся частей, и так далее. Если таким образом можно реализовать исходное разбиение, то разрезание будет тривиальным.

Обозначим через~\fo{A(x_1, \dots, x_n)} множество всех таких чисел~$z$, что прямоугольник с отношением сторон~$z$ можно тривиально разрезать на прямоугольники с отношениями сторон~$x_1$,~\dots,~$x_n$. Аналогично через~\fo{B(x_1, \dots, x_n)} обозначим множество всех таких чисел~$z$, что прямоугольник с отношением сторон~$z$ можно разрезать на прямоугольники с отношениями сторон~$x_1$,~\dots,~$x_n$ (не обязательно тривиальным образом). Очевидно, что~\fo{A(x_1, \dots, x_n)\subset B(x_1, \dots, x_n)}.

\glava\textbf{ Доказательство основной теоремы для тривиальных разрезаний.}\end{center}

Основную теорему мы докажем сначала для случая тривиальных разрезаний. Это означает, что в формулировке теоремы нужно заменить слово <<разрезать>> на словосочетание <<тривиально разрезать>>. Доказательство для этого случая является более наглядным, чем для произвольных разбиений. Далее в статье мы докажем эту теорему и для любых разрезаний.

Итак, прежде всего нам понадобятся несколько вспомогательных утверждений, сформулированных в виде лемм 1--6.

\lemma \textit{Множество \fo{A(x_1, \dots, x_n)} замкнуто относительно сложения, то есть для любых чисел \fo{a,b\in A(x_1, \dots, x_n)} выполнено \fo{(a+b)\in A(x_1, \dots, x_n)}.}\label{l1}

$\blacktriangleleft\,$Действительно, пусть \fo{a,b\in A(x_1, \dots, x_n)}. Заметим, что это означает, что из прямоугольников с отношениями сторон~$x_1$,~\dots,~$x_n$ мы можем тривиально составить любые прямоугольники с отношениями сторон $a$~и~$b$, в том числе прямоугольник со сторонами 1~и~$a$ и прямоугольник со сторонами 1~и~$b$. Приложим эти два прямоугольника друг к другу по стороне~1. Получим прямоугольник со сторонами 1~и~$a+b$.$\,\blacktriangleright$

\lemma\label{l2} \textit{Множество~\fo{A(x_1, \dots, x_n)} замкнуто относительно операции взятия обратного числа по умножению, то есть для любого числа~\fo{a\in A(x_1, \dots, x_n)} выполнено~\fo{a^{-1}\in A(x_1, \dots, x_n)}.}

$\blacktriangleleft\,$Действительно, по определению, если число $a$ --- отношение сторон некоторого прямоугольника, то число $a^{-1}$ также является отношением сторон того же прямоугольника.$\,\blacktriangleright$

\lemma\label{l3} \textit{Множество~\fo{A(x_1, \dots, x_n)} замкнуто относительно умножения на положительные рациональные числа, то есть для любого числа~\fo{a\in A(x_1, \dots, x_n)} и для любого числа~\fo{q\in\q^+} выполнено~\fo{aq\in A(x_1, \dots, x_n)}.}

\begin{wrapfigure}[9]{l}{0.26\linewidth} 
\vspace{-6.5ex}
\begin{picture}(45,45) 
\put(5,33){\line(1,0){15}} 
\put(5,29){\line(1,0){15}}
\put(5,25){\line(1,0){15}}
\put(26,33){\line(1,0){15}} 
\put(26,29){\line(1,0){15}}
\put(26,25){\line(1,0){15}}
\put(5,5){\line(1,0){15}} 
\put(5,9){\line(1,0){15}}
\put(5,13){\line(1,0){15}}
\put(26,5){\line(1,0){15}} 
\put(26,9){\line(1,0){15}}
\put(26,13){\line(1,0){15}}
\put(5,33){\line(0,-1){10}}
\put(11,33){\line(0,-1){10}}
\put(17,33){\line(0,-1){10}}
\put(29,33){\line(0,-1){10}}
\put(35,33){\line(0,-1){10}}
\put(41,33){\line(0,-1){10}}
\put(5,15){\line(0,-1){10}}
\put(11,15){\line(0,-1){10}}
\put(17,15){\line(0,-1){10}}
\put(29,15){\line(0,-1){10}}
\put(35,15){\line(0,-1){10}}
\put(41,15){\line(0,-1){10}}
\put(21,33){\circle*{0.5}}
\put(23,33){\circle*{0.5}}
\put(25,33){\circle*{0.5}}
\put(21,5){\circle*{0.5}}
\put(23,5){\circle*{0.5}}
\put(25,5){\circle*{0.5}}
\put(5,17){\circle*{0.5}}
\put(5,19){\circle*{0.5}}
\put(5,21){\circle*{0.5}}
\put(41,17){\circle*{0.5}}
\put(41,19){\circle*{0.5}}
\put(41,21){\circle*{0.5}}
\put(3.5,4.5){$^1$}
\put(3.5,8.5){$^1$}
\put(3.5,24.5){$^1$}
\put(3.5,28.5){$^1$}
\put(7.5,32){$^a$}
\put(13.5,32){$^a$}
\put(37.5,32){$^a$}
\put(31.5,32){$^a$}
\put(20,2.1){$ma$}
\put(41.5,18.5){$n$}
\put(17,-4){\ris\label{ri2}}
\end{picture}
\end{wrapfigure}

$\blacktriangleleft\,$Действительно, пусть~\fo{a\in A(x_1, \dots, x_n)} и~${q=\frac{m}{n}}$, где~\fo{m,n\in\n}. Прямоугольник со сторонами~$n$ и~$ma$, отношение сторон которого равно~${\frac{ma}{n}=qa}$, можно тривиально разрезать на прямоугольники со сторонами~1 и~$a$ (см. рис.~\ref{ri2}), а у этих прямоугольников, в свою очередь, существует тривиальное разбиение на прямоугольники с отношениями сторон~$x_1$,~\dots,~$x_n$, так как~\fo{a\in A(x_1, \dots, x_n)}.$\,\blacktriangleright$
 
\lemma\label{l4} \textit{Пусть \fo{(a-b\sqrt{p})\in\hor{p}\cap\R^+}. Тогда из прямоугольников с отношением сторон \fo{x=a+b\sqrt{p}} можно тривиально сложить любой прямоугольник с рациональным положительным отношением сторон и прямоугольник с отношением сторон \fo{a-b\sqrt{p}}.}

$\blacktriangleleft\,$Сначала докажем, что можно сложить прямоугольник с отношением сторон \fo{a-b\sqrt{p}}. По лемме~\ref{l2} имеем ${\frac{a-b\sqrt{p}}{a^2-pb^2}=\frac{1}{a+b\sqrt{p}}=\frac{1}{x}\in A(x)}$. Так как \fo{(a^2-pb^2)\in\q^+}, то по лемме~\ref{l3} получаем ${(a-b\sqrt{p})\in A(x)}$, что и требовалось. 

Теперь докажем, что можно сложить любой прямоугольник с рациональным положительным отношением сторон. По лемме~\ref{l1} прямоугольник с отношением сторон~$2a$ можно тривиально сложить из двух прямоугольников с отношениями сторон~\fo{a-b\sqrt{p}} и~\fo{a+b\sqrt{p}} (так как~\fo{(a-b\sqrt{p})+(a+b\sqrt{p})=2a}). Из прямоугольников с отношением сторон~\fo{2a\in \q^+} по лемме~\ref{l3} можно тривиальнно сложить прямоугольник с любым рациональным положительным отношением сторон.$\,\blacktriangleright$

\lemma\label{l5} \textit{Пусть \fo{(a-b\sqrt{p})\in\hor{p}\cap\R^-}. Тогда из прямоугольников с отношением сторон \fo{x=a+b\sqrt{p}} можно тривиально сложить любой прямоугольник с отношением сторон \fo{q\sqrt{p}} (где \fo{q\in\q^+}) и прямоугольник с отношением сторон \fo{b\sqrt{p}-a}.}

Доказательство леммы \ref{l5} аналогично доказательству леммы \ref{l4} --- оставляем его читателю.

\lemma\label{l6} \textit{Пусть множество~$P$ замкнуто относительно операции сложения и операции взятия обратного по умножению. Тогда если~\fo{x_1, \dots, x_n\in P}, то~\fo{P\supset A(x_1, \dots, x_n)}.}

$\blacktriangleleft\,$Прямоугольники с отношениями сторон~$x_1$,~\dots,~$x_n$ будем называть \textit{базовыми}. Доказательство провед\1м индукцией по количеству базовых прямоугольников в разбиении.

База индукции. Из одного базового прямоугольника можно тривиально сложить лишь прямоугольники с отношениями сторон~${x_1, \dots, x_n,\frac{1}{x_1}, \dots, \frac{1}{x_n}}$. По условию леммы~\fo{x_1, \dots, x_n\in P}. По другому условию леммы множество~$P$ замкнуто относительно операции взятия обратного по умножению --- поэтому~${\frac{1}{x_1}, \dots, \frac{1}{x_n}\in P}$. Значит, база индукции выполняется.

Шаг индукции. Предположим, все прямоугольники, которые можно тривиально сложить не более чем из~$k$ прямоугольников с отношениями сторон~$x_1$,~\dots,~$x_n$, имеют отношения сторон, принадлежащие множеству~\fo{P}. Любой прямоугольник~$ABCD$, тривиально составленный из~${k+1}$ базового прямоугольника, можно разрезать на два прямоугольника, тривиально составленные из базовых прямоугольников. Каждый из этих двух прямоугольников будет
\begin{wrapfigure}[9]{l}{0.28\linewidth} 
\vspace{-4.4ex} 
\begin{picture}(45,45) 
\put(5,10){\line(0,1){25}} 
\put(5,35){\line(1,0){40}}
\put(20,10){\line(0,1){25}}
\put(5,10){\line(1,0){40}}
\put(45,10){\line(0,1){25}}
\put(2,6.5){$A$}
\put(2,35.5){$B$}
\put(45,6.5){$D$}
\put(45,35.5){$C$}
\put(18.5,6.5){$F$}
\put(18.5,35.5){$E$}
\put(18,1){\ris\label{ri3}}
\end{picture}
\end{wrapfigure}
разрезан не более чем на~$k$ базовых прямоугольников, значит, по предположению индукции отношения сторон каждого из них принадлежат множеству~\fo{P}. Найд\1м отношения сторон прямоугольника~$ABCD$ (см. рис.~\ref{ri3}). Он разрезан на два прямоугольника ---~$ABEF$ и~$ECDF$. Так как~${\frac{BE}{AB},\frac{EC}{CD}\in P}$, то~${\frac{BC}{AB}=\frac{BE}{AB}+\frac{EC}{CD}\in P}$ (из замкнутости $P$ по сложению) и~${\frac{AB}{BC}=\left( \frac{BC}{AB}\right)^{-1}\in P}$ (из замкнутости~\fo{P} относительно операции взятия обратного по умножению). 

Таким образом, мы доказали методом математической индукции, что~\fo{P\supset A(x_1, \dots, x_n)}.$\,\blacktriangleright$

\textbf{Доказательство основной теоремы для тривиальных разрезаний.} Докажем первый пункт теоремы. По условию этого пункта существуют два таких числа $i$ и $j$, что~\fo{1\leqslant i,j\leqslant n},~\fo{a_i-b_i\sqrt{p}>0}, а~\fo{c_j-d_j\sqrt{p}<0}. 

Докажем, что~\fo{P:=\hor{p}\cap\R^+\subset A(x_1, \dots, x_n)}. Возьмём произвольное~${z=(e+f\sqrt{p})\in P}$. Рассмотрим~3 случая.

Случай 1:~\fo{e>0},~\fo{f>0}. По леммам~\ref{l4} и~\ref{l5} имеем~\fo{e,f\sqrt{p}\in A(x_1, \dots, x_n)}, следовательно,~\fo{e+f\sqrt{p}\in A(x_1, \dots, x_n)} по лемме~\ref{l1}.

Случай 2:~\fo{e<0} (очевидно, что при этом~\fo{f>0} и~\fo{pf^2-e^2>0}). Тогда~${\frac{1}{z}=\frac{-e+f\sqrt{p}}{pf^2-e^2}\in\horp{p}\subset A(x_1, \dots, x_n)}$ (по случаю~1), и по лемме~\ref{l2} получаем, что~\fo{z\in A(x_1, \dots, x_n)}.

Случай 3:~\fo{f<0} (очевидно, что при этом~\fo{e>0} и~\fo{e^2-pf^2>0}). Тогда~${\frac{1}{z}=\frac{e-f\sqrt{p}}{e^2-pf^2}\in\horp{p}\subset A(x_1, \dots, x_n)}$ (по случаю~1), и по лемме~\ref{l2} получаем, что~\fo{z\in A(x_1, \dots, x_n)}.

Очевидно, что множество~\fo{P} замкнуто относительно операций сложения и взятия обратного по умножению. Тогда, так как~\fo{x_1, \dots, x_n\in P}, то по лемме~\ref{l6} выполнено~\fo{P\supset A(x_1, \dots, x_n)}.

Мы доказали, что~\fo{P\supset A(x_1, \dots, x_n)} и~\fo{P\subset A(x_1, \dots, x_n)}, следовательно,~\fo{P=A(x_1, \dots, x_n)}. Первый пункт теоремы доказан.

Докажем второй пункт теоремы. Покажем сначала, что~\fo{M(x_1, \dots, x_n)\subset A(x_1, \dots, x_n)}. Для определ\1нности будем считать, что~${|b_1|/a_1=\underset{i}{\max}(|b_i|/a_i)}$,~\fo{1\leqslant i\leqslant n}. Если~\fo{b_1=0}, то все числа~$x_1$,~\dots,~$x_n$ являются рациональными,~\fo{M(x_1, \dots, x_n)=\q^+}, а~\fo{\q^+\subset A(x_1, \dots, x_n)} по лемме~\ref{l3}. В дальнейшем будем считать, что~\fo{b_1\neq 0}.

Пусть~\fo{z=(e+f\sqrt{p})\in M(x_1, \dots, x_n)}. По условию теоремы~${e\in\q^+}$, ${f\in\q}$ и~${|f|/e\leqslant |b_1|/a_1}$. Рассмотрим два случая.

Случай 1: ${f\in\q^+_0}$. Докажем, что $$z=e+f\sqrt{p}=\frac{f}{|b_1|}(a_1+|b_1|\sqrt{p})+\left( e-\frac{fa_1}{|b_1|}\right)\in A(x_1, \dots, x_n).$$

Действительно,~${\left( e-\frac{fa_1}{|b_1|}\right)\in\q^+_0}$ (так как~${\frac{|f|}{e}\leqslant\frac{|b_1|}{a_1}}$) и~${\frac{f}{|b_1|}\in\q^+_0}$ (так как~${f\in\q^+_0}$). По лемме~\ref{l4} выполнено~${(a_1-b_1\sqrt{p})\in A(x_1, \dots, x_n)}$, следовательно,~${(a_1+|b_1|\sqrt{p})\in A(x_1, \dots, x_n)}$. Также по лемме~\ref{l4} для любых~\fo{q\in\q^+} выполнено~\fo{q\in A(x_1, \dots, x_n)}. Применяя леммы~\ref{l1} и~\ref{l3}, получаем требуемое.

Случай 2:~${f\in\q^-}$. Докажем, что $$z=e+f\sqrt{p}=\frac{|f|}{|b_1|}(a_1-|b_1|\sqrt{p})+\left( e-\frac{|f|a_1}{|b_1|}\right)\in A(x_1, \dots, x_n).$$

Действительно,~${\left( e-\frac{|f|a_1}{|b_1|}\right)\in\q^+_0}$ (так как ${\frac{|f|}{e}\leqslant\frac{|b_1|}{a_1}}$) и ${\frac{|f|}{|b_1|}\in\q^+}$. Применяя леммы~\ref{l1} и~\ref{l3}, получаем требуемое.

Таким образом,~\fo{z\in A(x_1, \dots, x_n)}, то есть~\fo{M(x_1, \dots, x_n)\subset A(x_1, \dots, x_n)}.

Теперь проверим включение~\fo{M(x_1, \dots, x_n)\supset A(x_1, \dots, x_n)}. По лемме~\ref{l6} для этого достаточно проверить замкнутость множества~$M(x_1, \dots, x_n)$ относительно операций сложения и взятия обратного по умножению.

Докажем сначала, что~${M(x_1, \dots, x_n)}$ замкнуто по сложению. Действительно, пусть~\fo{(e+f\sqrt{p}),(g+h\sqrt{p})\in M(x_1, \dots, x_n)}. Тогда~${|f|\leqslant e\frac{|b|}{a}}$ и~${|h|\leqslant g\frac{|b|}{a}}$. Сложив эти неравенства, получим~${|f+h|\leqslant |f|+|h|\leqslant (e+g)\frac{|b|}{a}}$. Кроме того, очевидно,~\fo{(e+g)\in\q^+} и~\fo{(f+h)\in\q}. Следовательно,~\fo{[(e+g)+(f+h)\sqrt{p}]\in M(x_1, \dots, x_n)}.

Теперь покажем, что~$M(x_1, \dots, x_n)$ замкнуто относительно операции взятия обратного по умножению. Действительно, пусть~\fo{(e+f\sqrt{p})\in M(x_1, \dots, x_n)}. Тогда~${\frac{1}{e+f\sqrt{p}}=\frac{e-f\sqrt{p}}{e^2-pf^2}}$. Очевидно, что это число также принадлежит множеству~\fo{M(x_1, \dots, x_n)}.

Таким образом,~\fo{M(x_1, \dots, x_n)\supset A(x_1, \dots, x_n)} и~\fo{M(x_1, \dots, x_n)\subset A(x_1, \dots, x_n)}, следовательно,~\fo{M(x_1, \dots, x_n)=A(x_1, \dots, x_n)}. Второй пункт нашей теоремы доказан.

Доказательство третьего пункта теоремы почти дословно повторяет доказательство второго пункта (вместо леммы~\ref{l4} нужно использовать лемму~\ref{l5}).

\glava\textbf{ Доказательство основной теоремы.}\end{center}

Приступим к доказательству основной теоремы в общем случае, то есть без предположения тривиальности разрезаний. Для этого нам снова потребуется доказать две вспомогательные леммы и дать одно определение.

\lemma\label{l7} \textit{Если~\fo{x_1, \dots, x_n\in\hor{p}\cap\R^+}, то~\fo{B(x_1, \dots, x_n)\subset\hor{p}\cap\R^+}.} 

$\blacktriangleleft\,$Действительно, если~\fo{x_1, \dots, x_n\in\hor{p}} и \fo{z\in B(x_1, \dots, x_n)}, то по теореме~\ref{teorD} получаем, что число $z$ можно выразить через числа \fo{x_1, x_2, \ldots, x_n} с помощью сложения, вычитания, умножения и деления. Поскольку множество~\fo{\hor{p}} замкнуто относительно этих операций, то~\fo{z\in\hor{p}}. Число~$z$ --- отношение сторон прямоугольника, поэтому по определению~\fo{z>0}, следовательно, \fo{B(x_1, \dots, x_n)\subset\hor{p}\cap\R^+}.$\,\blacktriangleright$

Поэтому основную теорему достаточно доказать при условии, когда \fo{z\in\hor{p}\cap\R^+}, потому что это условие необходимо выполнено как при \fo{z\in M(x_1, \dots, x_n)} или \fo{z\in N(x_1, \dots, x_n)}, так и при \fo{z\in B(x_1, \dots, x_n)} (по лемме~\ref{l7}). Далее будем считать, что \fo{e+f\sqrt{p}=z\in\hor{p}\cap\R^+}.

\opred Пусть \fo{A,B,C\in\R}. \textit{<<Площадью>>} $S$ прямоугольника со сторонами ${\alpha+\beta\sqrt{p}}$ и ${\gamma+\delta\sqrt{p}}$, где \fo{\alpha,\beta,\gamma,\delta,p\in\q} и \fo{\sqrt{p}\notin\q}, называется число $$S:=\alpha\gamma A+\beta\gamma B+\alpha\delta B+\beta\delta C.$$

Легко проверить, что данное определение \2площади\3 не зависит от перестановки сторон прямоугольника.

\lemma\label{l9} \textit{Если прямоугольник разрезан на прямоугольники со сторонами, принадлежащими множеству \fo{\hor{p}}, то \2площадь\3 разрезаемого прямоугольника равна сумме \2площадей\3 прямоугольников, на которые он разрезан.}

\begin{wrapfigure}[6]{l}{0.25\linewidth} 
\vspace{-1ex}
\begin{picture}(45,83)
\thicklines
\put(0,80){\line(1,0){45}}
\put(14,72){\line(1,0){31}}
\put(21,76){\line(1,0){24}}
\put(0,60){\line(1,0){45}}
\put(0,60){\line(0,1){20}}
\put(14,60){\line(0,1){20}}
\put(21,72){\line(0,1){8}}
\put(25,60){\line(0,1){12}}
\put(35,60){\line(0,1){12}}
\put(45,60){\line(0,1){20}}
\thinlines
\multiput(0,72)(2,0){8}{\line(1,0){1}}
\multiput(0,76)(2,0){11}{\line(1,0){1}}
\multiput(21,60)(0,2){6}{\line(0,1){1}}
\multiput(25,80)(0,-2){5}{\line(0,-1){1}}
\multiput(35,80)(0,-2){5}{\line(0,-1){1}}
\put(15.5,55){\ris\label{rrr5}}
\end{picture}
\end{wrapfigure}
$\blacktriangleleft\,$Нетрудно убедиться, что сумма <<площадей>> двух прямоугольников со сторонами, принадлежащими множеству \fo{\hor{p}} и имеющих общую сторону, равна <<площади>> их объединения.

В общем случае сделаем продолжение каждого разреза, как показано на рис.~\ref{rrr5}. Очевидно, что тогда каждый прямоугольник нового разрезания будет также иметь квадратичные иррациональные стороны. Применяя утверждение предыдущего абзаца, получаем требуемое.$\,\blacktriangleright$

\textbf{Доказательство основной теоремы в случае \fo{{\bf \textit{n}=1}}.} Возможность разрезаний, когда \fo{a_1-b_1\sqrt{p}>0} и \fo{z\in M(x_1)} или когда \fo{a_1-b_1\sqrt{p}<0} и \fo{z\in N(x_1)}, нами уже доказана: по основной теореме для тривиальных разрезаний, когда выполнена одна из этих двух систем условий, существует тривиальное разрезание прямоугольника с отношением сторон $z$ на прямоугольники с отношением сторон $x_1$.

Остаётся доказать, что разрезание невозможно ни когда \fo{a_1-b_1\sqrt{p}>0} и \fo{z\notin M(x_1)}, ни когда \fo{a_1-b_1\sqrt{p}<0} и \fo{z\notin N(x_1)}. Предположим противное --- пусть прямоугольник с отношением сторон $z$ можно разрезать на прямоугольники с отношением сторон $x_1$ и выполнена одна из двух систем условий: либо \fo{a_1-b_1\sqrt{p}>0} и \fo{|f|/e>|b_1|/a_1}, либо \fo{a_1-b_1\sqrt{p}<0} и \fo{|e|/f>|a_1|/b_1}.

Будем называть разрезаемый прямоугольник с отношением сторон $z$ \textit{большим}, а прямоугольники с отношением сторон $x_1$, на которые он разрезается --- \textit{маленькими}.

Пусть большой прямоугольник имеет стороны 1 и \fo{e+f\sqrt{p}}. По теореме о силах тока и длинах сторон и лемме о сопротивлении цепи из~\cite{sobe_kva} стороны маленьких прямоугольников в этом случае должны выражаться через число $x_1$ с помощью сложения, вычитания, умножения и деления. Тогда, поскольку множество \fo{\hor{p}} замкнуто относительно этих операций, то стороны маленьких прямоугольников являются квадратичными иррациональностями. Таким образом, мы можем воспользоваться введённым выше определением \2площади\3 для большого и маленьких прямоугольников.

Рассмотрим \2площадь\3, для которой $$A=f,\,\,\,B=-e,\,\,\,C=\frac{2fa_1^2}{b_1^2}-pf.$$ \2Площадь\3 большого прямоугольника в этом случае равна $$S_{\mbox{\textit{\scriptsize{Б}}}}=e\cdot1\cdot f+f\cdot1\cdot(-e)+e\cdot0\cdot(-e)+(-e)\cdot0\cdot\left(\frac{2fa_1^2}{b_1^2}-pf\right)=0.$$

Теперь оценим \2площадь\3 маленьких прямоугольников. Пусть стороны какого-то маленького прямоугольника равны ${\alpha+\beta\sqrt{p}}$ и ${\gamma+\delta\sqrt{p}}$ и \fo{\frac{\gamma+\delta\sqrt{p}}{\alpha+\beta\sqrt{p}}=a_1+b_1\sqrt{p}}. Тогда \fo{\gamma+\delta\sqrt{p}=\alpha a_1+p\beta b_1+\sqrt{p}(\beta a_1+\alpha b_1)} и выполнены равенства $$\gamma=\alpha a_1+p\beta b_1\,\,\,\mbox{ и }\,\,\,\delta=\beta a_1+\alpha b_1.$$ Поэтому \2площадь\3 маленького прямоугольника записывается в виде $$S_{\mbox{\textit{\scriptsize{М}}}}=\alpha\gamma f-\beta\gamma e-\alpha\delta e+\beta\delta\left(\frac{2fa_1^2}{b_1^2}-pf\right)=$$
$$=\alpha f(\alpha a_1+p\beta b_1)-\beta e(\alpha a_1+p\beta b_1)-\alpha e(\beta a_1+\alpha b_1)+\beta\left(\frac{2fa_1^2}{b_1^2}-pf\right)(\beta a_1+\alpha b_1)=$$
$$=\alpha^2(fa_1-eb_1)+2\alpha\beta(\frac{fa_1^2}{b_1}-ea_1)+\beta^2(\frac{2fa_1^3}{b_1^2}-pfa_1-peb_1).$$

Заметим, что по предположению \fo{fa_1-eb_1\neq 0}. Поэтому выражение \fo{S_{\mbox{\textit{\scriptsize{М}}}}/\beta^2} --- квадратный трёхчлен относильно \fo{\alpha/\beta}. Покажем, что его дискриминант $D$ отрицателен. Тем самым будет доказано, что при фиксированных числах $a_1$, $b_1$, $e$, $f$ величина \fo{S_{\mbox{\textit{\scriptsize{М}}}}} либо всегда положительна, либо всегда отрицательна, а значит, сумма всех \2площадей\3 маленьких прямоугольников не равна нулю, то есть не равна площади большого прямоугольника. Таким образом, мы получим противоречие.

Итак, $$\frac{D}{4}=\frac{f^2a_1^4}{b_1^2}-\frac{2fea_1^3}{b_1}+e^2a_1^2-\frac{2f^2a_1^4}{b_1^2}+pf^2a_1^2+pfa_1eb_1+\frac{2fea_1^3}{b_1}-pfa_1eb_1-pe^2b_1^2=$$
$$=\frac{-f^2a_1^4}{b_1^2}+e^2a_1^2+pf^2a_1^2-pe^2b_1^2=(a_1^2-pb_1^2)\left( e^2-\frac{f^2a_1^2}{b_1^2}\right).$$ Имеем два случая: 1)~\fo{a_1-b_1\sqrt{p}<0} и \fo{|e|/f>|a_1|/b_1\Leftrightarrow |e|>f|a_1|/b_1}; 2)~\fo{a_1-b_1\sqrt{p}>0} и \fo{|f|/e>|b_1|/a_1\Leftrightarrow e<a_1|f|/|b_1|}. В обоих этих случаях \fo{D<0}.$\,\blacktriangleright$

\textbf{Доказательство основной теоремы.} Докажем первый пункт теоремы. По основной теореме для тривиальных разрезаний \fo{B(x_1, \dots, x_n)\supset A(x_1, \dots, x_n)=\hor{p}\cap\R^+}. Также по лемме~\ref{l7} выполняется обратное включение. Таким образом, \fo{B(x_1, \dots, x_n)=\hor{p}\cap\R^+}, что и требовалось в первом пункте доказываемой теоремы.

Докажем второй пункт теоремы. Не теряя общности, можем предположить, что~${|b_1|/a_1=\underset{i}{\max}(|b_i|/a_i)}$,~\fo{1\leqslant i\leqslant n}. Докажем, что тогда \fo{B(x_1)=B(x_1, \dots, x_n)}. По основной теореме для тривиальных разрезаний \fo{x_2, \dots, x_n\in A(x_1)}. Пусть \fo{z\in B(x_1, \dots, x_n)}. В разрезании прямоугольника с отношением сторон~$z$ на прямоугольники с отношениями сторон~$x_1$,~\dots,~$x_n$ заменим каждый прямоугольник с отношением сторон, принадлежащим множеству \fo{\lbrace x_2, \dots, x_n\rbrace}, на его тривиальное разрезание на прямоугольники с отношением сторон~$x_1$. Получим разрезание прямоугольника с отношением сторон~$z$ на прямоугольники с отношением сторон~$x_1$. Это мы доказали, что \fo{B(x_1, \dots, x_n)\subset B(x_1)}. Обратное включение очевидно следует из определения множества \fo{B(x_1, \dots, x_n)}. Таким образом, \fo{B(x_1)=B(x_1, \dots, x_n)}. По основной теореме для случая \fo{n=1} $${\displaystyle B(x_1, \dots, x_n)=B(x_1)=\left\lbrace  e+f\sqrt{p}\,\vert\, e\in\q^+, f\in\q,\frac{|f|}{e}\leqslant\frac{|b_1|}{a_1} \right\rbrace.}$$ Вспоминая, что мы предположили, что ~${|b_1|/a_1=\underset{i}{\max}(|b_i|/a_i)}$,~\fo{1\leqslant i\leqslant n}, получаем $${\displaystyle B(x_1, \dots, x_n)=\left\lbrace  e+f\sqrt{p}\,\vert\, e\in\q^+, f\in\q,\frac{|f|}{e}\leqslant\underset{i}{\max}\frac{|b_i|}{a_i} \right\rbrace,}$$ что и требовалось.

Третий пункт теоремы доказывается аналогично второму.

\glava\textbf{ Заключение.}\end{center}

Основная теорема является обобщением теоремы~\ref{TeoFLR}, но на первый взгляд они выглядят совершенно разными. Покажем, что на самом деле при \fo{n=1} эти теоремы эквивалентны.

\textbf{Утверждение.} \textit{Пусть \fo{u=\alpha+\beta\sqrt{p}} и \fo{v=\delta u+\gamma} (где \fo{\alpha,\beta,\gamma,\delta,p\in\q} и \fo{\beta\sqrt{p}\notin\q}) --- положительные квадратичные иррациональности. Тогда одна из двух систем условий $$\mbox{1}^\circ\mbox{) }\alpha-\beta\sqrt{p}>0\mbox{ и }v\in M(u)\,\,\,\mbox{ или}\,\,\,\mbox{ 2}^\circ\mbox{) }\alpha-\beta\sqrt{p}<0\mbox{ и }v\in N(u)$$ выполнена тогда и только тогда, когда выполнена одна из двух систем условий $$\mbox{1*) }\gamma=0\mbox{ и }\delta>0\,\,\,\mbox{ или}\,\,\,\mbox{ 2*) }\alpha\neq 0,\mbox{ }\frac{\gamma(\alpha^2-\beta^2p)}{\alpha}>0\mbox{ и }\delta+\frac{\gamma}{2\alpha}\geqslant 0.$$}

$\blacktriangleleft\,$Рассмотрим два случая: когда \fo{\alpha-\beta\sqrt{p}>0} и когда \fo{\alpha-\beta\sqrt{p}<0}.

1) Пусть
\begin{equation}
\label{us1}
\alpha-\beta\sqrt{p}>0.
\end{equation}
Сразу заметим, что это условие равносильно выполнению условия \fo{\alpha>|\beta|\sqrt{p}}, так как ${\alpha+\beta\sqrt{p}=u>0}$. Далее в этом пункте доказательства мы произведём ряд равносильных преобразований, которые верны лишь при условии~(\ref{us1}).

Итак, условие \fo{v\in M(u)}, очевидно, равносильно следующим двум условиям:
\begin{equation}
\label{f1}
\delta\alpha+\gamma>0 \mbox{  и  } \displaystyle{\frac{|\beta\delta|}{\delta\alpha+\gamma}\leqslant\frac{|\beta|}{\alpha}}.
\end{equation}

Теперь рассмотрим два подслучая: 1.1) \fo{\gamma=0} и 1.2) \fo{\gamma\neq 0}.

1.1) Если \fo{\gamma=0}, то выполнение системы условий (\ref{f1}) равносильно выполнению условия \fo{\delta>0}.

1.2) Если \fo{\gamma\neq 0}, то, пользуясь первым условием из (\ref{f1}), второе можно преобразовать так:
$$\displaystyle{\frac{|\beta\delta|}{\delta\alpha+\gamma}\leqslant\frac{|\beta|}{\alpha}}
\Leftrightarrow
\displaystyle{|\delta|\leqslant\delta+\frac{\gamma}{\alpha}}
\Leftrightarrow
\left\{
   \begin{array}{lcl}
      \displaystyle{-2\delta\leqslant\frac{\gamma}{\alpha}}\\
      \displaystyle{0\leqslant\frac{\gamma}{\alpha}}
   \end{array}
\right.
\Leftrightarrow
\left\{
   \begin{array}{lcl}
      \displaystyle{\delta+\frac{\gamma}{2\alpha}\geqslant 0}\\
      \gamma>0
   \end{array}
\right.
\Leftrightarrow
\left\{
\begin{array}{lcl}
         \alpha\neq 0\\
         \displaystyle{\delta+\frac{\gamma}{2\alpha}\geqslant 0}\\
         \displaystyle{\frac{\gamma(\alpha^2-\beta^2p)}{\alpha}>0,}
      \end{array}
      \right.$$
где на последнем шаге мы воспользовались предположением \fo{\alpha>|\beta|\sqrt{p}}.

Получаем, что в случае (1) система условий~\fo{(1^\circ)} выполнена тогда и только тогда, когда выполнена одна из двух систем условий~(1*) или~(2*).

2) Пусть
\begin{equation}
\label{us2}
\alpha-\beta\sqrt{p}<0.
\end{equation}
Это условие равносильно выполнению условия \fo{\alpha<|\beta|\sqrt{p}}, так как ${\alpha+\beta\sqrt{p}=u>0}$. Далее в этом пункте доказательства мы произведём ряд равносильных преобразований, которые верны лишь при условии~(\ref{us2}).

Итак, условие \fo{v\in N(u)}, очевидно, равносильно следующим двум условиям:
\begin{equation}
\label{f2}
\delta\beta>0 \mbox{  и  } \displaystyle{\frac{|\delta\alpha+\gamma|}{\beta\delta}\leqslant\frac{|\alpha|}{\beta}.}
\end{equation}

Снова рассмотрим два подслучая: 2.1) \fo{\gamma=0} и 2.2) \fo{\gamma\neq 0}.

2.1) Если \fo{\gamma=0}, то выполнение системы условий~(\ref{f2}) равносильно выполнению условия \fo{\delta>0}.

2.2) Если же \fo{\gamma\neq 0}, то, пользуясь первым условием из~(\ref{f2}), второе можно преобразовать так:
$$\displaystyle{\frac{|\delta\alpha+\gamma|}{\beta\delta}\leqslant\frac{|\alpha|}{\beta}}
\Leftrightarrow
      |\delta\alpha+\gamma|\leqslant\delta|\alpha|
\Leftrightarrow
         \left\{
         \begin{array}{lcl}
            \alpha\neq 0\\
            \displaystyle{|\delta+\frac{\gamma}{\alpha}|\leqslant\delta}
         \end{array}
         \right.
\Leftrightarrow
         \left\{
         \begin{array}{lcl}
            \alpha\neq 0\\
            \displaystyle{\delta+\frac{\gamma}{2\alpha}\geqslant 0}\\
            \displaystyle{\frac{\gamma}{\alpha}\leqslant 0}
         \end{array}
         \right.
\Leftrightarrow
\left\{
      \begin{array}{lcl}
         \alpha\neq 0\\
         \displaystyle{\delta+\frac{\gamma}{2\alpha}\geqslant 0}\\
         \displaystyle{\frac{\gamma(\alpha^2-\beta^2p)}{\alpha}>0.}
      \end{array}
      \right.$$
      
Получаем, что в случае (2) система условий~\fo{(2^\circ)} выполнена тогда и только тогда, когда выполнена одна из двух систем условий~(1*) или~(2*).

\glava\textbf{ Что дальше?}\end{center}

Естественным образом у читателя могут возникнуть несколько вопросов, касающихся обобщения полученного нами результата. Ниже приведены два наиболее интересных вопроса.

1. Какие прямоугольники можно разрезать на прямоугольники, подобные $n$ данным, если отношения сторон данных прямоугольников --- кубические иррациональности?

2. Какие прямоугольники можно разрезать на прямоугольники, подобные $n$ данным, если отношения сторон данных прямоугольников --- квадратичные иррациональности ${x_1=a_1+b_1\sqrt{p_1}}$,~\dots,~${x_n=a_n+b_n\sqrt{p_n}}$?

Автор не видит простого способа обобщить основную теорему и на эти случаи. Проблема заключается в нахождении аналогов лемм~\ref{l4} и~\ref{l5} в этих случаях. Эти вопросы пока ещё не решены.



В заключение приведём ещё несколько близких открытых вопросов.

1. Какие прямоугольники можно разрезать на подобные $n$ данным?

2. Какие многоугольники можно разрезать на прямоугольники с данным отношением сторон~$x$?

3. Какие многоугольники можно разрезать на трапеции, гомотетичные данной?\\

Автор благодарен М. Б. Скопенкову за постоянное внимание к данной работе.

\bibliographystyle{elsarticle-num}